\documentclass[12pt]{article}
\usepackage{mathrsfs}
\usepackage{epic,eepic,epsf,epsfig}
\usepackage{amsfonts,srcltx,mathrsfs}
\usepackage{amsmath}
\usepackage{amssymb}
\usepackage{amsbsy}
\usepackage{float}
\usepackage{graphicx}
\usepackage{amsfonts}
\usepackage{color}

\usepackage{url}
\newtheorem{lem}{Lemma}[section]
\newtheorem{theorem}[lem]{Theorem}
\newtheorem{defi}[lem]{Definition}

\newtheorem{case}{Case}
\newtheorem{prop}[lem]{Proposition}

\newtheorem{conjecture}[lem]{Conjecture}

\def\r{\rho} \def\s{\sigma}

 \def\lg{\langle} \def\rg{\rangle}
\def\nd{\mathrel{\bigm|\kern-.7em/}}

\def\f{\noindent}

\def\PSL{\hbox{\rm PSL}}

 \def\G{\hbox{\rm G}}

\def\Aut{\hbox{\rm Aut}}

\def\Aut{\hbox{\rm Aut}}

\def\demo{\f {\bf Proof.}\hskip10pt}

\newcommand{\qed}{\mbox{\raisebox{0.7ex}{\fbox{}}} \vspace{4truemm}}
\def\mz{{\mathbb Z}}

\begin{document}

\title{Existence of regular $3$-hypertopes with $2^n$ chambers}

\author{ \\ Dong-Dong Hou$^{a}$, Yan-Quan Feng*$^{a}$, Dimitri Leemans$^{b}$\\
$^{a}${\small Department of Mathematics, Beijing Jiaotong University, Beijing,
100044, P.R. China}\\
$^{b}${\small {D\'epartement de Math\'ematique, Universit\'e Libre de Bruxelles, 1050 Bruxelles Belgium}}}

\date{}
\maketitle

\footnotetext{*Corresponding author.\\
\mbox{}\hskip 0.7cm E-mails: holderhandsome$@$bjtu.edu.cn, yqfeng$@$bjtu.edu.cn, dleemans$@$ulb.ac.be}

\begin{abstract}
{For any positive integers $n, s, t, l$ such that $n \geq 10$, $s, t \geq 2$, $l \geq 1$ and $n \geq s+t+l$, a new  infinite family of regular 3-hypertopes
with type $(2^s, 2^t, 2^l)$ and automorphism group of order $2^n$ is constructed.
}

\bigskip
\f {\bf Keywords:} Regular hypertope, $2$-group, automorphism group.
\end{abstract}

\section{Introduction}
Polyhedra and their generalisation to higher ranks, polytopes, have been studied for more than two millenniums. Abstract regular polytopes generalize the concept of "realisable" polytopes, making them combinatorial objects consisting of a poset satisfying a series of axioms that the usual polytopes did.

Abstract polytopes, when regular, are in one-to-one correspondence with string C-groups, that are smooth quotients of Coxeter groups~\cite{ARP} . These groups have a Coxeter diagram that is a string. Dropping the string condition on the Coxeter diagram associated to a string C-group, Fernandes, Leemans and Weiss introduced hypertopes in~\cite{FDA2016} as a generalisation of abstract polytopes. Hypertopes are thin residually connected incidence geometries and if they are regular, one can associate to them a C-group (not necessary string). However, given a C-group, it is not always the case that the incidence geometry associated to it (via Tits' algorithm) is a hypertope.

Hypertopes of rank three, also called hyperhedra, are non-degenerate hypermaps. By that, we mean they satisfy the diamond condition, that is, given two incident hyperfaces $F$ and $G$, there are exactly two hyperfaces incident to both $F$ and $G$.

Hypertopes and C-groups are becoming an active topic currently. Ens~\cite{EE2018} classified regular toroidal hypertopes of rank four, and it was shown that their
automorphism groups are the quotients of infinite irreducible Coxeter groups of euclidean type with 4 generators. Catalano et al.~~\cite{CFHL} constructed an infinite family of hypertopes of rank four having the complete graph $K_4$, and their groups of rotational symmetries are isomorphic to the simple group $\PSL(2,q)$ with $q=p$ or $p^2$ where $p$ is a prime number satisfying some extra conditions.
Fernandes and Leemans~\cite{FL2018} classified C-groups of
rank $n-1$ and $n-2$ for the symmetric group $S_n$, and all these C-groups correspond to regular hypertopes. Fernandes et al.~\cite{FDA2018} investigated finite rank 4 structures obtained by hexagonal extensions of toroidal hypermaps, and many new examples that are regular or chiral were given even when the extensions are polytopal, which contains a new infinite family of finite nonlinear hexagonal extensions of the tetrahedron.

However, there are few results for nilpotent groups. If a regular hypertope has a nilpotent automorphism group, then it must be a 2-group. In this paper we focus on regular hypertopes of rank three whose automorphism groups are 2-groups. We prove the existence of regular 3-hypertopes
with type $(2^s, 2^t, 2^l)$ and automorphism group of order $2^n$ for $n \geq 10$, $s, t \geq 2$, $l \geq 1$ and $n\geq s+t+l$. In fact, we construct a group $G$ generated by three involutions $\r_0, \r_1, \r_2$, such that $G$ has order $2^n$, $\r_0\r_1, \r_1\r_2, \r_0\r_2$ have orders $2^s, 2^t, 2^l$ respectively, the pair $(G, \{\r_0, \r_1, \r_2\})$ satisfies the intersection property (see~2.1 for definition) and the subgroups $G_0=\lg \r_1, \r_2\rg$, $G_1=\lg \r_0, \r_2\rg$ and $G_2=\lg \r_0, \r_1\rg$
satisfy one of the Tits conditions (see Propostion~\ref{Tits}).

The paper is organized as follows. In Section~\ref{backgroud}, we give background definitions and properties of regular hypertopes and their automorphism groups needed to understand in this paper.
The main result and its proof is given in Section~\ref{Main Results}.

\section{Background results}\label{backgroud}

\subsection{C-groups}

Let $G$ be a group and $S=\{\rho_0,\cdots,\rho_{d-1}\}$ be a generating set of involutions of $G$.
For $I\subseteq \{0, \ldots, d-1\}$, let $G_I := \langle \rho_j : j \in I\rangle$.
We say that the pair $(G,S)$ satisfies the {\em intersection property} if $G_I \cap G_J = G_{I\cap J}$ for every $I,J\subseteq \{0, \ldots, d-1\}$. Clearly, the intersection property implies that $S$ is a minimal generating set of $G$.

A {\em C-group} is a pair $(G,S)$ satisfying the intersection property, where $G$ is a group and $S$ is a generating set of involutions of $G$, and the {\em rank} of $(G,S)$ is the cardinality of $S$.
A C-group is a {\em string C-group} if its set of generators $S$ can be ordered in such a way that $S:=\{\rho_0,\ldots,\rho_{d-1}\}$ satisfies the string property: $(\r_i\r_j)^2=1$ for all $i, j \in \{0, \cdots, d-1\}$ with $|i-j|>1$.  If $(G,S)$  only satisfies the string property, it is called a {\em string group generated by involutions} or \emph {sggi}.
It is known that string $C$-groups are in one-to-one correspondence with abstract regular polytopes~\cite[Section 2E]{ARP}.

The following proposition is called the {\em quotient criterion} for a string C-group.

\begin{prop}{\rm \cite[Section 2E]{ARP}}\label{stringC}
Let $(\G,\{\r_0, \r_1, \r_2\})$ be an sggi, and let $\Lambda = (\lg \s_{0}, \s_{1}, \s_{2}\rg,$ $\{\s_0,\s_1,\s_2\})$  be a string C-group.
If the mapping $\r_j \mapsto \s_j$ for $j=0 ,1, 2$ induces a homomorphism $\pi : \G \rightarrow \Lambda$, which is one-to-one on the subgroup
 $\lg \r_0, \r_1\rg$ or on $\lg \r_1, \r_2\rg$, then $(\G, \{\r_0, \r_1, \r_2\})$ is also a string $C$-group.
\end{prop}

The following proposition gives some string C-groups with type $\{4, 4\}$. It is proved in \cite[Section 8.3]{HW} and \cite{ARP}.

\begin{prop}\label{type44}
For $b \geq 2$, let
\begin{itemize}\setlength{\parskip}{-3pt}
  \item [$M_1$]$=\lg \r_0, \r_1, \r_2 \ |\ \r_0^2, \r_1^2, \r_2^2, (\r_0\r_1)^{4}, (\r_1\r_2)^{4}, (\r_0\r_2)^2, (\r_2\r_1\r_0)^{2b}\rg$,
  \item [$M_2$]$=\lg \r_0, \r_1, \r_2 \ |\ \r_0^2, \r_1^2, \r_2^2, (\r_0\r_1)^{4}, (\r_1\r_2)^{4}, (\r_0\r_2)^2, (\r_1\r_2\r_1\r_0)^b\rg$.
\end{itemize}
Then $|M_1|=16b^2$ and $|M_2|=8b^2$.
\end{prop}
In Proposition~\ref{type44}, it is easy to see that $o(\r_2\r_1\r_0)=2b$ in $M_1$ and $o(\r_1\r_2\r_1\r_0)=b$ in $M_2$.

\subsection{Regular hypertopes}

Hypertopes are a natural generalization of polytopes. In order to give the definition, we start with the definition of an {\em incidence system}~\cite{BC2013}.
\begin{defi}
{\rm
An {\em incidence system} $\Gamma :=(X, *, t, I)$ is a 4-tuple such that
\begin{itemize}\setlength{\parskip}{-3pt}
\item $X$ is a set whose elements are called the {\em elements} of $\Gamma$;
\item $I$ is a set whose elements are called the {\em types} of $\Gamma$;
\item $t: X \rightarrow I$ is a {\em type function}, associating to each elements $x \in X$ of $\Gamma$ a {\em type} $t(x) \in I$;
\item $*$ is a binary relation on $X$ called {\em incidence}, that is reflexive, symmetric and such that for all $x, y \in X$, if $x*y$ and
$t(x)=t(y)$ then $x=y$.
\end{itemize}
}
\end{defi}

The {\em incidence graph} of $\Gamma$ is the graph whose vertex set is $X$ and where two vertices are joined provided the corresponding elements of $\Gamma$ are incident. A {\em flag} $F$ is a set of pairwise incident elements of $\Gamma$, i.e. a clique of its incidence graph and the {\em type} of $F$
is $\{t(x) : x \in F\}$. A {\em chamber} is a flag of type $I$. An element $x$ is said to be {\em incident} to a flag $F$ when $x$ is incident to all elements of $F$
and we write $x*F$. An incidence system $\Gamma$ is a {\em incidence geometry} provided that every flag of $\Gamma$ is contained in a chamber. The {\em rank} of $\Gamma$ is the number of types of $I$, namely the cardinality of $I$.

Let $\Gamma :=(X, *, t, I)$ be an incidence geometry and $F$ a flag of $\Gamma$. The {\em residue} of $F$ in $\Gamma$ is the incidence geometry
$\Gamma_{F} :=(X_{F}, *_{F}, t_{F}, I_{F})$ where
\begin{itemize}\setlength{\parskip}{-3pt}
\item $X_{F}:=\{x \in X: x*F, x\notin F\}$;
\item $I_{F}=I \backslash t(F)$;
\item $t_{F}$ and $*_{F}$ are the restrictions of $t$ and $*$ to $X_{F}$ and $I_{F}$.
\end{itemize}

If each residue of rank at least $2$ of $\Gamma$ has a connected incidence graph then $\Gamma$ is said to be {\em residually connected}.
Moreover, $\Gamma$
is {\em thin} when every residue of rank one of $\Gamma$ contains exactly two elements. A {\em hypertope} is a residually connected thin incidence geometry.

Let $\Gamma :=(X, *, t, I)$ be an incidence system. An {\em automorphism} of $\Gamma$ is a mapping
$\alpha : (X, I) \rightarrow (X, I): (x, t(x)) \mapsto (\alpha(x), t(\alpha(x)))$ where
\begin{itemize}\setlength{\parskip}{-3pt}
\item $\alpha$ is a bijection on $X$;
\item for each $x, y \in X$, $x*y$ if and only if $\alpha(x)*\alpha(y)$;
\item for each $x, y \in X$, $t(x)=t(y)$ if and only if $t(\alpha(x))=t(\alpha(y))$.
\end{itemize}

An automorphism $\alpha$ of $\Gamma$ is called {\em type preserving} when for each $x \in X$, $t(\alpha(x))=t(x)$.
The set of type-preserving automorphisms of $\Gamma$ forms a group, denoted by $\Aut_{I}(\Gamma)$, and the set of automorphisms of $\Gamma$ also forms a group, denoted by $\Aut(\Gamma)$.

An incidence geometry $\Gamma$ is {\em flag-transitive} if $\Aut_{I}(\Gamma)$ is transitive on all flags of a given type $J$
for each type $J \subseteq I$. An incidence geometry $\Gamma$ is {\em chamber-transitive} if $\Aut_{I}(\Gamma)$ is transitive on all chambers of $\Gamma$.
Moreover, an incidence geometry $\Gamma$ is {\em regular} if $\Aut_{I}(\Gamma)$ acts regularly on the chambers (i.e. the action is semi-regular and transitive). The following proposition is folklore in incidence geometry.

\begin{prop}{\rm ~\cite[Proposition 2.2]{FDA2016}}\label{folklore}
Let $\Gamma$ be an incidence geometry. Then $\Gamma$ is chamber-transitive if and only if $\Gamma$ is flag-transitive.
\end{prop}

By Proposition~\ref{folklore}, a regular hypertope is a flag-transitive hypertope.
A rank one hypertope is a geometry with two elements. The polygons are precisely the hypertopes of rank two.
In fact, every abstract regular
polytope is a regular hypertope, and in rank three and higher, there are regular hypertopes that are not abstract regular polytopes.

Given an incidence system $\Gamma$ and a chamber $C$ of $\Gamma$, we may associate to the pair $(\Gamma, C)$ a pair consisting of a group $G$ and
a set $\{G_i : i \in I\}$ of subgroups of $G$ where $G :=\Aut_{I}(\Gamma)$ and $G_i$ is the stabilizer in $G$ of the element of type $i$ in $C$.
The following proposition shows how to reverse this construction, that is starting from a group and some of its subgroups, how to construct an incidence system.

\begin{prop}{\rm (Tits, 1956)~\cite{Tits1957}} \label{CoG}
Let $n$ be a positive integer and $I:=\{1, \cdots, n\}$. Let $G$ be a group together with a family of subgroups $(G_i)_{i \in I}$, $X$ the set consisting of all cosets $G_{i}g$ with $g \in G$ and $i \in I$, and $t: X\rightarrow I$ defined by $t(G_ig)=i$. Define an incidence relation $*$ on $X \times X$ by:
\begin{center}
$G_ig_1*G_jg_2$ iff $G_ig_1 \cap G_jg_2 \neq \emptyset$.
\end{center}
Then the $4$-tuple $\Gamma :=(X, *, t, I)$ is an incidence system having a chamber. Moreover, the group $G$ acts by right multiplication as an automorphism group on $\Gamma$. Finally, the group $G$ is transitive on the flags of rank less than $3$.
\end{prop}

For a group $G$ and a family of subgroups $(G_i)_{i \in I}$, the geometry $\Gamma$ constructed in Proposition~\ref{CoG} is called a {\em coset geometry}.
The following result gives a way to check whether a coset geometry (and
in particular a hypertope) is flag-transitive. See also Dehon~\cite{MD1994} or \cite{FL2018}.

\begin{prop}{\rm (Buekenhout, Hermand)~\cite{BH1991}}\label{FT}
Let $\mathcal{P}(I)$ be the set of all the subsets of $I$ and let $\alpha : \mathcal{P}(I) \backslash \{\emptyset\} \rightarrow I$ be a function such that
$\alpha(J) \in J$ for every $J \subset I, J \ne \emptyset$. Then $G$ is flag-transitive on coset geometry $\Gamma(G; (G_i)_{i \in I})$ if and only if,
for every $J \subset I$ such that $|J| \geq 3$, we have
$$\bigcap \limits_{j \in J-\alpha(J)} (G_jG_{\alpha(J)})=(\bigcap\limits_{j \in J-\alpha(J)} G_j)G_{\alpha(J)}.$$
\end{prop}

When $\Gamma(G; (G_i)_{i \in I})$ is a regular hypertope, from~\cite[Lemma 3.4]{FDA2016} we know that $\cap_{j \in I \backslash \{i\}}G_{j}$ is a group of order $2$, and define $\r_{i}$ as the unique involution in $\cap_{j \in I \backslash \{i\}}G_{j}$. We call the set $\{\r_i : i \in I\}$ the {\em distinguished generators} of $\Gamma(G; (G_i)_{i \in I})$.
\begin{prop}{\rm ~\cite[Theorem 4.1]{FDA2016}}
Let $I :=\{0, 1, \cdots, r-1\}$ and let $\Gamma(G; (G_i)_{i \in I})$ be a regular hypertope of rank $r$. Then the pair $(G, S)$ where $S$ is the set of distinguished generators of $\Gamma$ is a C-group of rank $r$.
\end{prop}

For now on, we construct a hypertope from a coset geometry $\Gamma(G; (G_i)_{i \in I})$ in a natural way, where $(G, S)$ is a C-group of rank $r$ with $I=\{0,1,\ldots, r-1\}$ and $G_i=G_{I \backslash \{i\}}$. For regular hypertopes, we have the following proposition.

\begin{prop}{\rm ~\cite[Theorem 4.6]{FDA2016} }\label{Hypertope}
Let $(G, \{\r_0, \r_1, \ldots, \r_{r-1}\})$ be a C-group of rank $r$, and let $G_i=\lg \r_j \ | \ \r_j \in S, j \in I \backslash \{i\} \rg$ for all $i \in I:= \{0, \ldots, r-1\}$. If $G$ is flag-transitive on $\Gamma(G; (G_i)_{i \in I})$, then $\Gamma$ is a regular hypertope.

\end{prop}

Let $(G, \{\r_0, \r_1, \r_2\})$ be a C-group of rank $3$ and let $\Gamma(G; (G_i)_{i \in I})$ be the coset geometry constructed from $(G, \{\r_0, \r_1, \r_2\})$, where $G_i=\lg \r_j \ | \ j \in \{0,1,2\}\backslash \{i\} \rg$. Let $I=\{0,1,2\}$ and let $\alpha : \mathcal{P}(I) \backslash \{\emptyset\} \rightarrow I$ be a function such that $\alpha (\{0, 1, 2\}) = \{0\}$, where $\mathcal{P}(I)$ is the set of all subsets of $I$. Clearly, $G_1G_0 \cap G_2G_0 = (G_1 \cap G_2)G_0$ if and only if $G_0G_1 \cap G_0G_2 = G_0(G_1 \cap G_2)$, and by Proposition~\ref{FT}, $\Gamma(G; (G_i)_{i \in I})$ is flag-transitive if and only if $G_1G_0 \cap G_2G_0 =(G_1 \cap G_2)G_0$, for which some equivalent conditions were given by Tits.

\begin{prop}{\rm (Tits)~\cite{Tits1974}}\label{Tits}
Let $G_0, G_1, G_2$ be three subgroups of a group $G$. Then the following conditions are equivalent.
\begin{itemize}\setlength{\parskip}{-3pt}
\item [{\rm(1)}] $G_0G_1 \cap G_0G_2 = G_0(G_1 \cap G_2)$;
\item [{\rm(2)}] $(G_0 \cap G_1)\cdot(G_0 \cap G_2)=(G_1G_2) \cap G_0$;
\item [{\rm(3)}] If the three cosets $G_0x, G_1y$ and $G_2z$ have pairwise nonempty intersection, then $G_0x \cap G_1y \cap G_2z\ne \emptyset$.
\end{itemize}

\end{prop}

\subsection{Group theory}

Let $G$ be a group. For $x,y\in G$, we use $[x,y]$ as an abbreviation for the
{\em commutator} $x^{-1}y^{-1}xy$ of $x$ and $y$, and $[H, K]$ for the subgroup generated by all commutators $[x, y]$ with $x \in H$ and $y \in K$, when $H$ and $K$ are subgroups of $G$.
The following  proposition is a basic property of commutators and its proof is straightforward.

\begin{prop}\label{commutator}
Let $G$ be a group. Then, for any $x, y, z \in G$, $[xy, z]=[x, z]^y[y, z]$ and $[x, yz]=[x, z][x, y]^z$.
\end{prop}

The {\em commutator (or derived)} subgroup $G'$ of a group $G$ is the subgroup generated by all commutators $[x, y]$ for any $x, y \in G$.
The {\em Frattini subgroup}, denoted by $\Phi(G)$, of a finite group $G$ is defined to be the intersection of all maximal subgroups of $G$. Let $G$ be a finite $p$-group for a prime $p$,  and set $\mho_1(G) = \lg g^p\ |\ g \in G\rg$. The following theorem is the well-known Burnside Basis Theorem.

\begin{prop}{\rm ~\cite[Theorem 1.12]{GroupBookss}}\label{burnside}
Let $G$ be a $p$-group and $|G: \Phi(G)| = p^d$.
\begin{itemize}
\item [(1)] $G/\Phi(G) \cong \mz_p^d$. Moreover, if $N \lhd G$ and $G/N$ is elementary abelian, then $\Phi(G) \leq N$.

\item [(2)] Every minimal generating set of $G$ contains exactly $d$ elements.

\item [(3)] $\Phi(G) = G' \mho_1(G)$. In particular, if $p=2$, then $\Phi(G) = \mho_1(G)$.
\end{itemize}
\end{prop}

By the above proposition, all minimal generating sets have the same cardinality, which is called the {\em rank} of $G$ and denoted by $d(G)$. This implies that a given $2$-group has only $C$-group representations with a fixed rank, that is, the rank of the $2$-group.

\section{Main Result}\label{Main Results}

To prove the main result, we need the following lemma.
\begin{lem}\label{group}
Let $M_1$ and $M_2$ be the groups defined in Proposition~\ref{type44}. Then
$M_1\cong (D_{2b} \times D_{2b}) \rtimes (C_2 \times C_2)$ and
$M_2\cong (D_{2b} \times D_{2b}) \rtimes C_2$.
\end{lem}

\demo Recall that
\begin{itemize}\setlength{\parskip}{-3pt}
  \item [$M_1$]$=\lg \r_0, \r_1, \r_2 \ |\ \r_0^2, \r_1^2, \r_2^2, (\r_0\r_1)^{4}, (\r_1\r_2)^{4}, (\r_0\r_2)^2, (\r_2\r_1\r_0)^{2b}\rg$,
  \item [$M_2$]$=\lg \r_0, \r_1, \r_2 \ |\ \r_0^2, \r_1^2, \r_2^2, (\r_0\r_1)^{4}, (\r_1\r_2)^{4}, (\r_0\r_2)^2, (\r_1\r_2\r_1\r_0)^b\rg$.
\end{itemize}

We first consider $M_1$. Write $A=\lg \r_1^{\r_0}, \r_1^{\r_2}\rg$ and $B=\lg \r_1, \r_1^{\r_0\r_2}\rg$. Since $(\r_2\r_1\r_0)^{2b}=1$ and $\r_0\r_2=\r_2\r_0$, we have $(\r_2\r_1\r_0\r_2\r_1\r_0)^b=
(\r_1^{\r_2}\r_1^{\r_0})^b$. By Proposition~\ref{type44}, $o(\r_2\r_1\r_0)=2b$, implying $o(\r_1^{\r_2}\r_1^{\r_0})=b$. It follows that $A=\lg \r_1^{\r_0}, \r_1^{\r_2}\rg \cong D_{2b}$. Since $\lg \r_1, \r_1^{\r_0\r_2}\rg =\lg \r_1^{\r_0}, \r_1^{\r_2}\rg^{\r_2}$, we have $B=\lg \r_1, \r_1^{\r_0\r_2}\rg\cong D_{2b}$. Note that $[\r_1, \r_1^{\r_0}]=(\r_1\r_0)^4=1$ and $[\r_1, \r_1^{\r_2}]=(\r_1\r_2)^4=1$. Furthermore, $[\r_1^{\r_0\r_2}, \r_1^{\r_0}]=[\r_1^{\r_2}, \r_1]^{\r_0}=1$ and $[\r_1^{\r_0\r_2}, \r_1^{\r_2}]=[\r_1^{\r_0}, \r_1]^{\r_2}=1$. It follows that $[A, B]=1$ and hence $AB\leq M_1$, that is, $AB$ is a subgroup of $M_1$.

Note that $B^{\r_0}=\lg \r_1, \r_1^{\r_0\r_2}\rg^{\r_0}=\lg \r_1^{\r_0}, \r_1^{\r_2}\rg=A$ and $A^{\r_0}=\lg \r_1^{\r_0}, \r_1^{\r_2}\rg^{\r_0}=\lg \r_1, \r_1^{\r_0\r_2}\rg=B$. Similarly, $A^{\r_2}=B$ and $B^{\r_2}=A$. It follows that $\lg\r_0,\r_2\rg$ normalizes $AB$ and hence $\lg\r_0,\r_2\rg AB\leq M_1$. Furthermore, $\lg\r_0,\r_2\rg AB=\lg\r_0,\r_2,\r_1^{\r_0}, \r_1^{\r_2}, \r_1, \r_1^{\r_0\r_2}\rg=\lg \r_0,\r_1,\r_2\rg=M_1$. By Proposition~\ref{type44}, $16b^2=|M_1|=|\lg\r_0,\r_2\rg AB|$, and since $\lg\r_0,\r_2\rg\cong C_2\times C_2$ and $|A|=|B|=4b^2$, we have $\lg\r_0,\r_2\rg\cap (AB)=1$ and $A\cap B=1$. Since $[A, B]=1$, we have $M_1=(A \times B) \rtimes \lg \r_0, \r_2\rg \cong (D_{2b} \times D_{2b}) \rtimes (C_2 \times C_2)$.

Now we consider $M_2$. Write $C=\lg \r_0, \r_2^{\r_1}\rg$ and $D=\lg \r_0^{\r_1}, \r_2\rg$. Since $(\r_1\r_2\r_1\r_0)^b=1$, we have $(\r_2^{\r_1}\r_0)^b=1$ and $(\r_2\r_0^{\r_1})^{b}=((\r_2^{\r_1}\r_0)^b)^{\r_1}=1$.
By Proposition~\ref{type44}, $o(\r_1\r_2\r_1\r_0)=b$, implying  $o(\r_2^{\r_1}\r_0)=o(\r_2\r_0^{\r_1})=b$.
It follows $C\cong D_{2b}$ and $D \cong D_{2b}$.
Note that $[\r_0, \r_0^{\r_1}]=(\r_0\r_1)^4=1, [\r_0^{\r_1}, \r_2^{\r_1}]=[\r_0, \r_2]^{\r_1}=1$ and $[\r_2^{\r_1}, \r_2]=(\r_1\r_2)^4=1$. Then $[C,D]=1$ and hence $CD \leq M_2$. Clearly, $C^{\r_1}=D$ and $D^{\r_1}=C$, implying $\lg \r_1 \rg CD\leq M_2$. It follows $\lg \r_1 \rg CD=\lg \r_1,\r_0,\r_2^{\r_1},\r_0^{\r_1},\r_2\rg=M_2$. Again by Proposition~\ref{type44}, $8b^2=|M_2|=|\lg \r_1 \rg CD|$, and hence $C\cap D=1$ and $\lg\r_1\rg\cap(CD)=1$. Since $[C,D]=1$, we have $M_2=(C \times D)\rtimes \lg \r_1\rg \cong (D_{2b} \times D_{2b}) \rtimes C_2$. \hfill\qed

After the paper was finished, we found that Lemma~\ref{group} was also proved  in \cite{ARP}, but with a different method. For completeness, we still keep the proof here.
The following is the main result of this paper.

\begin{theorem}\label{existmaintheorem}
For any positive integers $n, s, t, l$ such that $n \geq 10$, $s, t \geq 2, l \geq 1$ and $n \geq s+t+l$, let $R(\r_0, \r_1, \r_2)= \{\r_0^2, \r_1^2, \r_2^2, (\r_0\r_1)^{2^{s}}, (\r_1\r_2)^{2^{t}}, (\r_0\r_2)^{2^{l}}, [(\r_0\r_1)^4, \r_2], [\r_0,(\r_1\r_2)^4],$
$ [(\r_0\r_2)^2, \r_1]\}$ and let
\begin{small}
$$G=\left\{
\begin{array}{ll}
\lg \r_0, \r_1, \r_2 \ |\ R(\r_0, \r_1, \r_2), [(\r_0\r_1)^2, \r_2]^{2^{\frac{n-s-t-l}{2}}}\rg, & n-s-t-l\mbox{ even }\\
\lg \r_0, \r_1, \r_2 \ |\ R(\r_0, \r_1, \r_2), [(\r_0\r_1)^2, (\r_1\r_2)^2]^{2^{\frac{n-s-t-l-1}{2}}}\rg, & n-s-t-l\mbox{ odd. }
\end{array}
\right.$$
\end{small}
The group $G$ is the automorphism group of a regular 3-hypertope with type $(2^s, 2^t, 2^l)$ and order $2^n$. Moreover,
\begin{small}
$$G\cong\left\{
\begin{array}{ll}
(C_{2^{s-2}} \times C_{2^{t-2}} \times C_{2^{l-1}}).((D_{2^{\frac{n-s-t-l+3}{2}}} \times D_{2^{\frac{n-s-t-l+3}{2}}})\rtimes(C_2 \times C_2) ), & n-s-t-l\mbox{ even }\\

(C_{2^{s-2}} \times C_{2^{t-2}} \times C_{2^{l-1}}).((D_{2^{\frac{n-s-t-l+4}{2}}}\times D_{2^{\frac{n-s-t-l+4}{2}}})\rtimes C_2), & n-s-t-l\mbox{ odd. }
\end{array}
\right.$$
\end{small}
\end{theorem}

\demo Let $L_1=\lg \r_0, \r_1, \r_2 \ |\ \r_0^2, \r_1^2, \r_2,  (\r_0\r_1)^{2^s} \rg$, $L_2=\lg \r_0, \r_1, \r_2 \ |\  \r_0, \r_1^2, \r_2^2,  (\r_1\r_2)^{2^t} \rg$ and $L_3=\lg \r_0, \r_1, \r_2 \ |\ \r_0^2, \r_2^2, \r_1,  (\r_0\r_2)^{2^l} \rg$. Then $L_1\cong D_{2^{s+1}}$, $L_2\cong D_{2^{t+1}}$ and $L_3\cong D_{2^{l+1}}$.

Let $A=\lg (\r_0\r_1)^4\rg$, $B=\lg (\r_0\r_1)^4\rg$ and $C=\lg (\r_0\r_2)^2 \rg$.
Write $o(g)$ for the order of $g$ in $G$ for any $g\in G$.

\medskip
\f {\bf Claim 1}: In $G$, $o(\r_0\r_1)=2^s$, $o(\r_1\r_2)=2^t$, $o(\r_0\r_2)=2^l$ and $|G|=2^n$. Furthermore, $A\unlhd G$, $B\unlhd G$, $C\unlhd G$ and $ABC=A\times B\times C\cong C_{2^{s-2}}\times C_{2^{t-2}}\times C_{2^{l-1}}$.

Since $\langle\r_1,\r_2\rangle$ is a dihedral group, we have  $(\r_1\r_2)^{\r_1}=(\r_1\r_2)^{\r_2}=(\r_1\r_2)^{-1}$.
Since $[\r_0,(\r_1\r_2)^4]=1$ in $G$, we have $A\unlhd G$. Similarly, $B \unlhd G$ and $C \unlhd G$ as $[(\r_0\r_1)^4,\r_2]=1$ and $[(\r_0\r_2)^2, \r_1]=1$. It is easy to see that the generators $\r_0, \r_1, \r_2$ in $L_1$ satisfy all relations in $G$. Thus $L_1$ is an epimorphic image of $G$, and since $\r_0\r_1$ has order $2^s$ in $L_1$, we have $o(\r_0\r_1)=2^s$ in $G$. Similarly, $o(\r_1\r_2)=2^t$ and $o(\r_0\r_2)=2^l$ in $G$. It follows that $A\cong C_{2^{s-2}}$, $B\cong C_{2^{t-2}}$ and $C\cong C_{2^{l-1}}$.

Let $K=ABC$. Now we prove that $K=A\times B\times C$ and $|G/K|=2^{n-s-t-l+5}$. Write $R_1(\r_0, \r_1, \r_2)=\{\r_0^2, \r_1^2, \r_2^2, (\r_0\r_1)^{2^s}, (\r_1\r_2)^{2^t}, (\r_0\r_2)^{2}, [(\r_0\r_1)^4, \r_2], [\r_0, (\r_1\r_2)^4]\}$.
Since $C=\lg (\r_0\r_2)^2\rg\unlhd G$, we have $G/C\cong G_1$, where
$$G_1=\left\{
\begin{array}{ll}
\lg \r_0, \r_1, \r_2 \ |\  R_1(\r_0, \r_1, \r_2), [(\r_0\r_1)^2, \r_2]^{2^{\frac{n-s-t-l}{2}}}\rg, & n-s-t-l\mbox{ even }\\
\lg \r_0, \r_1, \r_2 \ |\ R_1(\r_0, \r_1, \r_2), [(\r_0\r_1)^2, (\r_1\r_2)^2]^{2^{\frac{n-s-t-l-1}{2}}}\rg, & n-s-t-l\mbox{ odd. }
\end{array}
\right.$$
Clearly, the generators $\r_0$, $\r_1$, $\r_2$ in $L_1$ satisfy all relations in $G_1$, and hence $\r_0\r_1$ has order $2^s$ in $G_1$. This implies $o(\r_0\r_1C)=2^s$ in $G/C$ and so $|AC/C|=o((\r_0\r_1C)^4)=2^{s-2}$. It follows that $|AC|=|C||AC/C|=2^{s+l-3}$, and since $|AC||A\cap C|=|A||C|=2^{s+l-3}$, we have $|A\cap C|=1$ and hence $AC=A\times C$. Similarly, the generators $\r_0$, $\r_1$, $\r_2$ in $L_2$ satisfy all relations in $G_1$, implying $o(\r_1\r_2C)=2^t$ in $G/C$.

Set $R_2(\r_0, \r_1, \r_2)= \{\r_0^2, \r_1^2, \r_2^2, (\r_0\r_1)^{4},
(\r_1\r_2)^{2^t}, (\r_0\r_2)^{2},  [\r_0,(\r_1\r_2)^4]\}$. Then $G/AC\cong G_2$, where
$$G_2=\left\{
\begin{array}{ll}
\lg \r_0, \r_1, \r_2 \ |\ R_2(\r_0, \r_1, \r_2), [(\r_0\r_1)^2, \r_2]^{2^{\frac{n-s-t-l}{2}}}\rg, & n-s-t-l\mbox{ even }\\
\lg \r_0, \r_1, \r_2 \ |\ R_2(\r_0, \r_1, \r_2), [(\r_0\r_1)^2, (\r_1\r_2)^2]^{2^{\frac{n-s-t-l-1}{2}}}\rg, & n-s-t-l\mbox{ odd. }
\end{array}
\right.$$
The generators $\r_0$, $\r_1$, $\r_2$ in $L_2$ satisfy all relations in $G_2$, and hence $\r_1\r_2$ has order $2^t$ in $G_2$, which implies $o(\r_1\r_2AC)=2^t$ and $|B(AC)/(AC)|=o((\r_1\r_2AC)^4)=2^{t-2}$. It follows that $|B(AC)|=|B(AC)/(AC)||AC|=2^{s+t+l-5}$, and since $|B(AC)||B\cap (AC)|=|B||AC|=2^{s+t+l-5}$, we have $|B\cap (AC)|=1$ and hence $K=ABC=A\times B\times C\cong  C_{2^{s-2}}\times C_{2^{t-2}}\times C_{2^{l-1}}$. In particular, $|K|=2^{s+t+l-5}$.

Set $R_3(\r_0, \r_1, \r_2)= \{\r_0^2, \r_1^2, \r_2^2, (\r_0\r_1)^{4}, (\r_1\r_2)^{4}, (\r_0\r_2)^{2} \}$. Then $G/K\cong G_3$, where
$$G_3=\left\{
\begin{array}{ll}
\lg \r_0, \r_1, \r_2 \ |\ R_3(\r_0, \r_1, \r_2), [(\r_0\r_1)^2, \r_2]^{2^{\frac{n-s-t-l}{2}}}\rg, & n-s-t-l\mbox{ even }\\
\lg \r_0, \r_1, \r_2 \ |\ R_3(\r_0, \r_1, \r_2), [(\r_0\r_1)^2, (\r_1\r_2)^2]^{2^{\frac{n-s-t-l-1}{2}}}\rg, & n-s-t-l\mbox{ odd. }
\end{array}
\right.$$

Assume that $n-s-t-l$ is even. Since $\r_0\r_2=\r_2\r_0$ and $(\r_0\r_1)^2=(\r_0\r_1)^{-2}=(\r_1\r_0)^2$, we have $[(\r_0\r_1)^2, \r_2]=(\r_1\r_0\r_1\r_2)^2$ in $G_3$, and so $[(\r_0\r_1)^2, \r_2]^{2^{\frac{n-s-t-l}{2}}}=(\r_1\r_0\r_1\r_2)^{2^{\frac{n-s-t-l+2}{2}}}$.
Note that $2^{\frac{n-s-t-l+2}{2}} \geq 2^{\frac{0+2}{2}}=2 $. By Proposition~\ref{type44}, $|G/K|=|G_3|=8\cdot (2^{\frac{n-s-t-l+2}{2}})^2=2^{n-s-t-l+5}$.

Assume that $n-s-t-l$ is odd. Noting that $(\r_0\r_1)^4=(\r_1\r_2)^4=1$ in $G_3$, we have
$1=[(\r_0\r_1)^2, (\r_1\r_2)^2]^{2^{\frac{n-s-t-l-1}{2}}}= (((\r_0\r_1)^2(\r_1\r_2)^2)^2)^{2^{\frac{n-s-t-l-1}{2}}}=
(((\r_0\r_1\r_2)^2)^2)^{2^{\frac{n-s-t-l-1}{2}}}$
$=(\r_0\r_1\r_2)^{2\cdot 2^{\frac{n-s-t-l+1}{2}}}$ because $\r_0\r_2=\r_2\r_0$.
Since $2^{\frac{n-s-t-l+1}{2}} \geq 2^{\frac{1+1}{2}}=2$,
Proposition~\ref{type44} implies $|G/K|=|G_3|=16 \cdot (2^{\frac{n-s-t-l+1}{2}})^2=2^{n-s-t-l+5}$.

In both cases, $|G/K|=|G_3|=2^{n-s-t-l+5}$ and hence  $|G|=|G/K|\cdot |K|=2^n$, as claimed. Since $K\cong  C_{2^{s-2}}\times C_{2^{t-2}}\times C_{2^{l-1}}$, Lemma~\ref{group} implies
\begin{small}
$$G\cong\left\{
\begin{array}{ll}
(C_{2^{s-2}} \times C_{2^{t-2}} \times C_{2^{l-1}}).((D_{2^{\frac{n-s-t-l+3}{2}}} \times D_{2^{\frac{n-s-t-l+3}{2}}})\rtimes(C_2 \times C_2) ), & n-s-t-l\mbox{ even }\\

(C_{2^{s-2}} \times C_{2^{t-2}} \times C_{2^{l-1}}).((D_{2^{\frac{n-s-t-l+4}{2}}}\times D_{2^{\frac{n-s-t-l+4}{2}}})\rtimes C_2), & n-s-t-l\mbox{ odd. }
\end{array}
\right.$$
\end{small}

\medskip
\f {\bf Claim 2:} $(G/C, \{\r_0C, \r_1C, \r_2C\})$ is a string C-group with $o(\r_0\r_2C)=2$, $o(\r_0\r_1C)=2^s$ and $o(\r_1\r_2C)=2^t$, and $(G, \{\r_0, \r_1, \r_2\})$ is a C-group.

Note that $K=ABC=A\times B\times C$ and $G/K\cong G_3$. By \cite[Section 8.3]{HW} and \cite{ARP}, $(G_3,\{\r_0,\r_1,\r_2\})$ and so $(G/K,\{\r_0K,\r_1K,\r_2K\})$ is a string C-group. Clearly, there is a natural epimorphism from $G/(AC)$ to $G/K$ induced by $\r_0AC\mapsto \r_0K$, $\r_1AC\mapsto \r_1K$ and $\r_2AC\mapsto \r_2K$ as $AC\leq K$. Note that $o(\r_0\r_1K)=4$ and $o(\r_0\r_2K)=2$ in $G/K$. Then $o(\r_0\r_1AC)=4$ and  $o(\r_0\r_2AC)=2$ in $G/(AC)$. Thus $G/AC$ is an sggi, and Proposition~\ref{stringC} implies that $(G/AC, \{\r_0AC, \r_1AC, \r_2AC\})$ is a string C-group. Recall that $o(\r_1\r_2AC)=2^t$ in  $G/AC$.

Similarly, there is a natural epimorphism from $G/C$ to $G/(AC)$ induced by  $\r_0C\mapsto \r_0AC$, $\r_1C\mapsto \r_1AC$ and $\r_2C\mapsto \r_2AC$, and hence $o(\r_1\r_2C)=2^t$ in $G/C$. Furthermore, $(G/C, \{\r_0C, \r_1C, \r_2C\})$ is a string C-group with $o(\r_0\r_2C)=2$ and $o(\r_0\r_1C)=2^s$ in $G/C$. In particular, $\{\r_0C, \r_1C, \r_2C\}$ is a minimal generating set of $G/C$ and by Proposition~\ref{burnside}, $G/C$ has rank $3$, that is, $d(G/C)=3$.

Write $G_0=\lg \r_1, \r_2\rg$, $G_1=\lg \r_0, \r_2\rg$ and $G_2=\lg \r_0, \r_1\rg$. To prove that $(G, \{\r_0, \r_1, \r_2\})$ is a C-group, we only need to show that $(G, \{\r_0, \r_1, \r_2\})$ satisfies the intersection property, which is equivalent to show that  $G_0 \cap G_2 = \lg \r_1 \rg$, $G_1 \cap G_2 = \lg \r_0 \rg$ and $G_0 \cap G_1 = \lg \r_2 \rg$ in $G$.

Suppose $G_0 \cap G_1 > \lg \r_2 \rg$. Since $G_0$ and $G_1$ are dihedral, we have
$(\r_1\r_2)^{2^{t-1}}\in G_0 \cap G_1$ and $(\r_0\r_2)^{2^{l-1}}\in G_0 \cap G_1$, both of which belong to the center of $G_0 \cap G_1$. Clearly, $G_0 \cap G_1$ has no subgroup isomorphic to $\mz_2^3$, implying that $\lg \r_1,(\r_1\r_2)^{2^{t-1}}\rg=\lg \r_1,(\r_0\r_2)^{2^{l-1}}\rg$.
It follows $(\r_1\r_2)^{2^{t-1}}=(\r_0\r_2)^{2^{l-1}}$ or $(\r_1\r_2)^{2^{t-1}}=\r_1(\r_0\r_2)^{2^{l-1}}$. For the former,  $o(\r_1\r_2C)\leq 2^{t-1}$ in $G/C$, contradicting $o(\r_1\r_2C)=2^t$. For the latter, $d(G/C)\leq 2$ as $G=\lg \r_1, \r_1\r_2, \r_0\r_2\rg$, contradicting $d(G/C)=3$.

Suppose  $G_1 \cap G_2 > \lg \r_0 \rg$. Similarly to the previous paragraph,
$(\r_0\r_2)^{2^{l-1}}=(\r_0\r_1)^{2^{s-1}}$ or $(\r_0\r_2)^{2^{l-1}}=\r_0(\r_0\r_1)^{2^{s-1}}$, which is impossible because the former implies $o(\r_0\r_1C)<2^s$ in $G/C$ and the latter implies $d(G/C)<3$.

Suppose $G_0 \cap G_2 > \lg \r_1 \rg$. Then
$(\r_0\r_1)^{2^{s-1}}=(\r_1\r_2)^{2^{t-1}}$ or $(\r_0\r_1)^{2^{s-1}}=\r_1(\r_1\r_2)^{2^{t-1}}$. The latter is impossible because $d(G)=3$. For the former, $(\r_0\r_1C)^{2^{s-1}}=(\r_1\r_2C)^{2^{t-1}}$
and hence $\lg \r_0C, \r_1C\rg\cap \lg \r_1C, \r_2C\rg>\lg \r_1C\rg$, contradicting that $G/C$ is a string C-group.

It follows that $G_0 \cap G_2 = \lg \r_1 \rg$, $G_1 \cap G_2 = \lg \r_0 \rg$ and $G_0 \cap G_1 = \lg \r_2 \rg$ in $G$, and hence $(G, \{\r_0, \r_1, \r_2\})$ is a C-group, as claimed.

\medskip
\f {\bf Claim 3}: $(G_0 \cap G_1)\cdot (G_0 \cap G_2)=(G_1 G_2) \cap G_0$.

Recall $G_0=\lg \r_1, \r_2\rg$, $G_1=\lg \r_0, \r_2\rg$ and $G_2=\lg \r_0, \r_1\rg$.
By Claim~2, $(G_0 \cap G_1)\cdot (G_0 \cap G_2)=\lg \r_2\rg \lg \r_1\rg=\{1, \r_2, \r_1, \r_2\r_1\}$.

Noting that $G_1$ and $G_2$ are dihedral, we have $G_1=\{(\r_2\r_0)^i\r_0^j\ |\ i\in \mathbb{Z}_{2^{l}}, j \in \mathbb{Z}_{2}\}$ and $G_2=\{\r_0^j(\r_0\r_1)^k\ |\ k\in \mathbb{Z}_{2^{s}}, j \in \mathbb{Z}_{2}\}$. It follows that
$G_1G_2=\{(\r_2\r_0)^i\r_0^j$
$(\r_0\r_1)^k \ | \ i\in \mathbb{Z}_{2^{l}}, j \in \mathbb{Z}_{2}, k \in \mathbb{Z}_{2^{s}} \}$. It is straightforward that $\{1, \r_2, \r_1, \r_2\r_1\}\subseteq (G_1G_2) \cap G_0$.

On the other hand, assume $x \in (G_1G_2) \cap G_0$. Then $x=(\r_2\r_0)^i\r_0^j(\r_0\r_1)^k=(\r_1\r_2)^{p}\r_1^{q}$ for some $ i\in \mathbb{Z}_{2^{l}}, k \in \mathbb{Z}_{2^{s}}, p \in \mathbb{Z}_{2^{t}}$ and $j, q \in \mathbb{Z}_{2}$.

Let $i$ be even. Then $xC=\r_0^j(\r_0\r_1)^kC=(\r_1\r_2)^{p}\r_1^{q}C$. By Claim~2, $G/C$ is a string C-group, and hence $xC\in \lg \r_0C, \r_1C\rg\cap \lg \r_1C, \r_2C\rg=\lg\r_1C\rg$. It follows $x\in C$ or $x\in \r_1C$. Note that $C=\lg (\r_0\r_2)^2\rg\leq G_1$. If $x\in C$ then $x\in C\cap G_0\leq G_1\cap G_0=\lg\r_2\rg$ because $G$ is a C-group (Claim~2), and since $\r_2\not\in C$, we have $x=1$. If $x\in \r_1C$ then $x\in \r_1G_1\cap G_0=\r_1(G_1\cap G_0)=\r_1\lg\r_2\rg$, and since $\r_1\r_2\not\in \r_1C$, we have $x=\r_1$.

Let $i$ be odd. Then $xC=\r_2\r_0^{j+1}(\r_0\r_1)^kC=(\r_1\r_2)^{p}\r_1^{q}C$, and hence $\r_2xC\in \lg \r_0C, \r_1C\rg\cap \lg \r_1C, \r_2C\rg=\lg\r_1C\rg$. By the previous paragraph, $\r_2x=1$ or $\r_1$, that is, $x=\r_2$ or $\r_2\r_1$.

It follows that $(G_1G_2) \cap G_0= \{1,\r_1,\r_2,\r_2\r_1\}$, and hence $(G_0 \cap G_1)\cdot (G_0 \cap G_2)=(G_1 G_2) \cap G_0$, as claimed.

\medskip
By Claim~3, $(G_0 \cap G_1)\cdot (G_0 \cap G_2)=(G_1 G_2) \cap G_0$, and by Proposition~\ref{Tits}, this is equivalent to having
$G_0G_1 \cap G_0G_2 = G_0(G_1 \cap G_2)$.
By Proposition~\ref{FT}, the coset geometry $\Gamma(G,\{G_0,G_1,G_2\})$ is therefore flag-transitive and by Proposition~\ref{Hypertope}, $\Gamma(G,\{G_0,G_1,G_2\})$ is a regular hypertope of type $(2^s,2^t,2^l)$ and with automorphism group $G$.
\hfill\qed

To end this paper, we would like to propose the following conjecture.

\begin{conjecture}
For any positive integers $n, s, t, l$ such that $n \geq 10$, $s, t \geq 2$, $l \geq 1$ and $n < s+t+l$, there is no regular 3-hypertope with type $(2^s, 2^t, 2^l)$ and  automorphism group of order $2^n$.
\end{conjecture}

By~\cite[Theorem 3.2]{SmallestPolytopes}, the conjecture is true for $l=1$, and {\sc Magma} shows that it is also true for $n\leq 10$. Furthermore, one may show that the conjecture is true if one of $\lg\r_0,\r_1\rg\lg\r_1,\r_2\rg$, $\lg\r_0,\r_1\rg\lg \r_0,\r_2\rg$ and $\lg\r_0,\r_2\rg\lg \r_1,\r_2\rg$ is a group.

\medskip

\f {\bf Acknowledgements:} This work was supported by the National Natural Science Foundation of China (11571035, 11731002) and the 111 Project of China (B16002).


\begin{thebibliography}{99}

\bibitem{GroupBookss} Y. Berkovich,
 Groups of Prime Power Order,
 vol. 1, Walter de Gruyter, Berlin, 2008.

\bibitem{BCP97}
W. Bosma, J. Cannon, C. Playoust,
 The {M}agma {A}lgebra {S}ystem. {I}:  the user language,
 J. Symbolic Comput. 24 (1997) 235--265.

\bibitem{BC2013}
F. Buekenhout, A. M. Cohen,
Diagram Geometry: Related to classical groups and buildings, A Series of Modern Surveys in Mathematics
[Results in Mathematics and Related Areas. 3rd Series. A Series of Modern Surveys in Mathematics],
Ergebnisse der Mathematik und ihrer Grenzglebiete. 3. Folge, vol. 57,
Springer, Heidelberg, 2013, pp. xiv+592.

\bibitem{BH1991}
F. Buekenhout, M. Hermand,
On flag-transitive geometries and groups,
Travaux de Math\'{e}matiques de l'Universit\'{e} Libre de Bruxelles 1 (1991) 45--78.

\bibitem{CFHL}
D. Catalano, M. E. Fernandes, I. Hubard, D. Leemans,
Hypertopes with tetrahedral diagram,
Electron. J. Combin. 25 (2018) \#P3.22.


\bibitem{SmallestPolytopes}
M. Conder,
  The smallest regular polytopes of given rank,
Adv. Math. 236 (2013) 92--110.


\bibitem{HW}
H. S. M. Coxeter, W. O. J. Moser,
Generators and Realitions for Discrete Groups, fourth ed.,
Springer-Verlag, New York, 1972.

\bibitem{MD1994}
M. Dehon,
Classifying geometries with {\sc Cayley},
J. Symbolic Comput. 17 (1994) 259--276.

\bibitem{EE2018}
E. Ens,
Rank 4 toroidal hypertopes,
Ars Math. Contemp. 15 (2018) 67--79.



\bibitem{FDA2016}
M. E. Fernandes, D. Leemans, 	A. I. Weiss,
Highly Symmetric Hypertopes,
Aequationes Math. 90 (2016) 1045--1067.

\bibitem{FDA2018}
M. E. Fernandes, D. Leemans, A. I. Weiss,
Hexagonal extensions of toroidal maps and hypermaps,
in M. D. E. et al. (Eds.), Discrete Geometry and Symmetry, Springer Proceedings in Mathematics and Statistics,
 2018, pp. 147--170.



\bibitem{FL2018}
M. E. Fernandes, D. Leemans,
C-groups of high rank for the symmetric groups,
J. Algebra 508 (2018) 196--218.

\bibitem{ARP}
P. McMullen, E. Schulte,
  Abstract regular polytopes,
 Cambridge University Press, Cambridge, 2002.

\bibitem{Tits1957}
J. Tits,
Sur les analogues alg\'{e}briques des groupes semi-simples complexes, in: Colloque d'alg\`ebre sup\'{e}rieure, tenu \`a Bruxelles du 19 au 22 d\'{e}cembre 1967,
Center Belgede Recherches Math\'{e}matiques, \'{E}tablissements Ceuterick,
Librairie Gauthier-Villars, Louvain, Paris, 1957, pp. 261-289.
\bibitem{Tits1974}
J. Tits,
Buildings of Spherical Type and Finite BN-pairs,
Springer, Berlin, 1974.

\end{thebibliography}
\end{document}